\documentclass{article}
\usepackage{amsmath,amssymb}
\begin{document}

\newtheorem{thm}{Theorem}[section]
\newtheorem{prop}[thm]{Proposition}
\newenvironment{dfn}{\medskip\refstepcounter{thm}
\noindent{\bf Definition \thesection.\arabic{thm}\ }}{\medskip}
\newenvironment{remark}[1][Remark]{\begin{trivlist}
\item[\hskip \labelsep {\bfseries #1}]}{\end{trivlist}}
\newenvironment{note}[1][Note]{\begin{trivlist}
\item[\hskip \labelsep {\bfseries #1}]}{\end{trivlist}}
\newenvironment{notes}[1][Notes]{\begin{trivlist}
\item[\hskip \labelsep {\bfseries #1}]}{\end{trivlist}}
\newenvironment{proof}[1][,]{\medskip\ifcat,#1
\noindent{{\it Proof}:\ }\else\noindent{\it Proof of #1.\ }\fi}
{\hfill$\square$\medskip}
\newenvironment{ack}[1][Acknowledgements]{\begin{trivlist}
\item[\hskip \labelsep {\bfseries #1}]}{\end{trivlist}}

\def\eq#1{{\rm(\ref{#1})}}
\def\Q{{\mathbb Q}}
\def\R{{\mathbb R}}
\def\O{{\mathbb O}}
\def\Z{{\mathbb Z}}
\def\GL{\mathop{\rm GL}}
\def\SL{\mathop{\rm SL}}
\def\SU{\mathop{\rm SU}}
\def\SO{\mathop{\rm SO}}
\def\U{\mathbin{\rm U}}
\def\sech{\mathop{\rm sech}\nolimits}
\def\sn{\mathop{\rm sn}\nolimits}
\def\cn{\mathop{\rm cn}\nolimits}
\def\dn{\mathop{\rm dn}\nolimits}
\def\d{{\rm d}}
\def\w{\wedge}
\def\C{{\mathbb C}}
\def\Re{\mathop{\rm Re}\nolimits}
\def\Im{\mathop{\rm Im}\nolimits}
\def\G2{\mathop\textrm{G}_2}
\def\Spin{\mathop{\rm Spin}}
\def\N{\mathbb{N}}

\title{Calibrated Submanifolds of $\R^7$ and $\R^8$ with Symmetries}
\author{\textsc{Jason Dean Lotay}\\ University College\\ Oxford}
\date{}

\maketitle

\section{Introduction}

In this article, we describe a method of constructing certain types of
calibrated submanifold of $\R^7$ and $\R^8$ with \emph{symmetries}.
The main result is the exhibition of explicit examples of
$\U(1)$-invariant \emph{associative} cones in $\R^7$ and
\emph{Cayley 4-folds} in $\R^8$ which are invariant under $\SU(2)$.
This research is motivated by the work of Joyce in \cite{Joy2} on
special Lagrangian (SL) $m$-folds in $\C^m$, and
the work of the author in \cite{Lotay1}.

In Section \ref{s2}, we describe the calibrations and calibrated
submanifolds that are the focus of our study.  These are called
\emph{associative 3-folds} and \emph{coassociative 4-folds} in
$\R^7$ and \emph{Cayley 4-folds} in $\R^8$.

The method of construction to produce calibrated submanifolds with
\emph{symmetries} is discussed in Section \ref{s3}.  The key result
is that we may define examples using a system of first-order
ordinary differential equations.  This section also reviews the
relevant material from \cite{Lotay1}.

Sections \ref{s4} and \ref{s5} contain the explicit examples.  The
first gives the system of differential equations defining
$\U(1)$-invariant associative cones.  These equations are solved in
a special case to give a 4-dimensional family of associative cones
over $T^2$.  Further, using the material in \cite[$\S$6]{Lotay1} and
these cones, we produce examples of \emph{ruled} associative
3-folds.

Section \ref{s5} considers Cayley 4-folds invariant under an action
of $\SU(2)$. The family of all Cayley 4-folds invariant under this
action is described using a real octic and three real quartics.  
Cayley 4-folds invariant under $\SU(2)$ are also considered in \cite{Gu}; there
is some overlap between our example and those given in this reference.

The final section gives some further examples of systems of ordinary
differential equations defining associative, coassociative and
Cayley submanifolds, each associated with a symmetry group which is
described.


\begin{notes}
\begin{itemize}\item[]
\item[(a)] Manifolds are assumed to be nonsingular and submanifolds
to be immersed unless stated otherwise.
\item[(b)] By a \emph{cone} in $\R^n$ we shall mean a dilation-invariant
submanifold of $\R^n$ which is nonsingular except possibly at $0$. 
\end{itemize}
\end{notes}

\section{Calibrated submanifolds of $\R^7$ and $\R^8$}\label{s2}

\subsection{Calibrated geometry}

We define \emph{calibrations} and \emph{calibrated submanifolds}
following the approach in \cite{HarLaw}.

\begin{dfn}\label{ch1s1dfn1}
Let $(M,g)$ be a Riemannian manifold.  An \emph{oriented tangent
$k$-plane} $V$ on $M$ is an oriented $k$-dimensional vector subspace
$V$ of $T_xM$, for some $x$ in $M$. Given an oriented tangent
$k$-plane $V$ on $M$, $g|_V$ is a Euclidean metric on $V$ and hence,
using $g|_V$ and the orientation on $V$, there is a natural volume
form, $\text{vol}_V$, which is a $k$-form on $V$.

A closed $k$-form $\eta$ on $M$ is a \emph{calibration} on $M$ if
$\eta | _V \leq\,\text{vol}_V$ for all oriented tangent $k$-planes
$V$ on $M$, where $\eta | _V = \kappa \cdot\text{vol}_V$ for some
$\kappa\in\R$, so $\eta | _V \leq\,\text{vol}_V$ if $\kappa\leq 1$.
An oriented $k$-dimensional submanifold $N$ of $M$ is a
\emph{calibrated submanifold} or \emph{$\eta $-submanifold} if $\eta
| _{T_xN} =\, \text{vol}_{T_xN}$ for all $x\in N$.
\end{dfn}

 Calibrated submanifolds are \textit{minimal} submanifolds
\cite[Theorem II.4.2]{HarLaw}.  The minimality of calibrated
submanifolds provides the following property, as discussed in
\cite{HarLaw}.

\begin{thm}\label{ch1s1thm1}
A calibrated submanifold is real analytic wherever it is
nonsingular.
\end{thm}

\subsection{Associative and coassociative submanifolds of $\R^7$}

The convention we adopt here for calibrations on $\R^7$ agree with
\cite[Chapter 10]{Joy1}.

\begin{dfn}\label{ch2s3dfn1} Let $(x_1,\ldots,x_7)$ be coordinates on $\R^7$ and write
$d{\bf x}_{ij\ldots k}$ for the form $dx_i\w dx_j\w\ldots\w dx_k$.
Define a 3-form $\varphi_0$ by:
\begin{equation}\label{ch2s3eq1}
\varphi_0 = d{\bf x}_{123}+d{\bf x}_{145}+d{\bf x}_{167}+d{\bf
x}_{246}- d{\bf x}_{257}-d{\bf x}_{347}-d{\bf x}_{356}.
\end{equation}
By \cite[Theorem IV.1.4]{HarLaw}, $\varphi_0$ is a calibration on
$\R^7$ and submanifolds calibrated with respect to $\varphi_0$ are
called \emph{associative 3-folds}.

The 4-form $\ast\varphi_0$, where $\varphi_0$ and $\ast\varphi_0$
are related by the Hodge star, is given by:
\begin{equation}\label{ch2s3eq2}
\ast\varphi_0 = d{\bf x}_{4567}+d{\bf x}_{2367}+d{\bf
x}_{2345}+d{\bf x}_{1357}-d{\bf x}_{1346}-d{\bf x}_{1256}-d{\bf
x}_{1247}.
\end{equation}
By \cite[Theorem IV.1.16]{HarLaw}, $\ast\varphi_0$ is a calibration
on $\R^7$, and $\ast\varphi_0$-submanifolds are called
\emph{coassociative 4-folds}.
\end{dfn}

\begin{remark} The form $\varphi_0$ is often referred to as the $\G2$ 3-form
on $\R^7$ since the Lie group $\G2$ may be defined as the stabilizer
of $\varphi_0$ in $\GL(7,\R)$.
\end{remark}

We have a far more useful description of coassociative 4-folds which
follows from \cite[Proposition IV.4.5 \& Theorem IV.4.6]{HarLaw}.

\begin{prop}\label{ch2s3prop3} A 4-dimensional submanifold $M$ of\/ $\R^7$,
with an appropriate orientation, is coassociative if and only if
$\varphi_0|_M\equiv 0$.
\end{prop}

\subsection{Cayley submanifolds of $\R^8$}

Our definition of a distinguished 4-form on $\R^8$ used to describe
Cayley 4-folds agrees with the convention in \cite[Chapter
10]{Joy1}.

\begin{dfn}\label{ch2s4dfn1} Let $(x_1,\ldots,x_8)$ be coordinates on $\R^8$  and write
$d{\bf x}_{ij\ldots k}$ for the form $dx_i\w dx_j\w\ldots\w dx_k$.
Define
 a 4-form $\Phi_0$ by:
\begin{align}
\Phi_0 &=  d{\bf x}_{1234}+d{\bf x}_{1256}+d{\bf x}_{1278}+d{\bf
x}_{1357}-d{\bf x}_{1368}-d{\bf x}_{1458}-d{\bf x}_{1467}
\nonumber \\
&{}+ d{\bf x}_{5678}+d{\bf x}_{3478}+d{\bf x}_{3456}+d{\bf
x}_{2468}-d{\bf x}_{2457}-d{\bf x}_{2367}-d{\bf
x}_{2358}.\label{ch2s4eq1}
\end{align}
\noindent By \cite[Theorem IV.1.24]{HarLaw}, $\Phi_0$ is a
calibration on $\R^8$, and submanifolds calibrated with respect to
$\Phi_0$ are called \emph{Cayley 4-folds}.
\end{dfn}

\begin{remark}
The stabilizer of $\Phi_0$ in $\GL(8,\R)$ is the Lie group
$\Spin(7)$.  We may thus refer to $\Phi_0$ as the $\Spin(7)$ 4-form.
\end{remark}

\section{Constructing examples with symmetries}\label{s3}

\subsection{Evolution equations}

In \cite{Lotay1}, an \emph{evolution equation} for associative
3-folds in $\R^7$ was derived as a generalisation of the work of
Joyce \cite{Joy2} on special Lagrangian $m$-folds in $\C^m$. The
proof relies on Theorem \ref{ch1s1thm1} and the following result
from Harvey and Lawson \cite[Theorem IV.4.1]{HarLaw}.

\begin{thm}\label{ch4s1thm2}
Let $P$ be a 2-dimensional real analytic submanifold of\/ $\R^7$.
There locally exists a real analytic associative 3-fold $N$ in
$\R^7$ which contains $P$. Moreover, $N$ is locally unique.
\end{thm}

We now present the theorem \cite[Theorem 4.3]{Lotay1}.

\begin{thm}\label{ch4s1thm3}
Let $P$ be a compact, orientable, 2-dimensional, real analytic
manifold, $\chi$ a real analytic nowhere vanishing section of\/
$\Lambda^2TP$ and $\psi:P \rightarrow\R^7$ a real analytic embedding
(immersion).  There exist $\epsilon > 0$ and a unique family
$\{\psi_t:t\in(-\epsilon,\epsilon)\}$ of real analytic maps
$\psi_t:P\rightarrow\R^7$ with $\psi_0  =  \psi$ satisfying
\begin{equation}\label{ch4s1eq1}
\left(\frac{d\psi_t}{dt}\right)^d   =
(\psi_t)_*(\chi)^{ab}(\varphi_0)_{abc}(g_0)^{cd}, 
\end{equation}

\noindent where $(g_0)^{cd}$ is the inverse of the Euclidean metric
on $\R^7$, using index notation for tensors on $\R^7$. Define
$\Psi:(-\epsilon, \epsilon) \times P \rightarrow\R^7$ by $\Psi(t,p)
= \psi_t(p)$. Then $M = \text{\emph{Image}}\,\Psi$ is a nonsingular
embedded (immersed) associative 3-fold in $\R^7$.
\end{thm}

We sketch the key ideas in the proof.  Since $P$ is compact and $P$,
$\chi$, $\psi$ are real analytic, the \emph{Cauchy--Kowalevsky
Theorem} \cite[Theorem B.1]{Racke} from the theory of partial
differential equations
 gives a family of maps $\psi_t$ as stated.  We may therefore define
 $\Psi$ and $M$ as in the statement of the theorem.
Theorem \ref{ch4s1thm2} implies there locally exists a locally
unique associative 3-fold $N$ containing $\psi(P)$. Showing that $N$
and $M$ agree near $\psi(P)$, using the fact that $\varphi_0$ is a
calibration, allows us to deduce that $M$ is associative.

\medskip

Using the associative case as a model we can quickly derive
analogous evolution equations for coassociative and Cayley 4-folds.

We first require two results, \cite[Theorem IV.4.3]{HarLaw} and
\cite[Theorem IV.4.6]{HarLaw}, which are both similar to Theorem
\ref{ch4s1thm2}.

\begin{thm}\label{ch5s1thm2}
Suppose $P$ is a 3-dimensional real analytic submanifold of\/$\R^7$
such that $\varphi_0 |_P \equiv 0$. There locally exists a real
analytic coassociative 4-fold $N$ in $\R^7$ which contains $P$.
Moreover, $N$ is locally unique.
\end{thm}

\begin{remark} Unlike Theorem \ref{ch4s1thm2}, we have to impose an extra condition
on the boundary submanifold $P$ in order to extend it to a
coassociative 4-fold in $\R^7$.\end{remark}

\begin{thm}\label{ch5s1thm3}
Suppose $P$ is a 3-dimensional real analytic submanifold of\/
$\R^8$. There locally exists a real analytic Cayley 4-fold $N$ in
$\R^8$ which contains $P$.  Moreover, $N$ is locally unique.
\end{thm}

With these results at our disposal, it is clear that we may prove
results like Theorem \ref{ch4s1thm3} for coassociative and Cayley
4-folds in exactly the same manner, so we omit the proofs.

\begin{thm}\label{ch5s1thm4}
Let $P$ be a compact, orientable, 3-dimensional, real analytic
manifold, $\chi$ a real analytic nowhere vanishing section of\/
$\Lambda^3 TP$ and $\psi: P\rightarrow\R^7$ a real analytic
embedding (immersion) such that $\psi^*(\varphi_0)\equiv 0$ on $P$.
There exist $\epsilon >0$ and a unique family $\{\psi_t:t\in
(-\epsilon,\epsilon)\}$ of real analytic maps
$\psi_t:P\rightarrow\R^7$ with $\psi_0=\psi$ satisfying
\begin{equation}\label{ch5s1eq1}
\left(\frac{d\psi_t}{dt}\right)^e =
(\psi_t)_*(\chi)^{abc}(\ast\varphi_0)_{abcd} (g_0)^{de}
\end{equation}
\noindent using index notation for tensors on $\R^7$, where
$(g_0)^{de}$ is the inverse of the Euclidean metric on $\R^7$.
Define $\Psi:(-\epsilon,\epsilon) \times P\rightarrow \R^7$ by
$\Psi(t,p)=\psi_t(p)$.  Then $M = \text{\emph{Image}}\, \Psi$ is a
nonsingular embedded (immersed) coassociative 4-fold in $\R^7$.
\end{thm}

\begin{note} The condition
$\psi^*(\varphi_0)|_P\equiv0$ implies that $\varphi_0$ vanishes on
the real analytic 3-fold $\psi(P)$ in $\R^7$ and allows us to apply
Theorem \ref{ch5s1thm2} as required.
\end{note}

\begin{thm}\label{ch5s1thm5}
Let $P$ be a compact, orientable, 3-dimensional, real analytic
manifold, $\chi$ a real analytic nowhere vanishing section of\/
$\Lambda^3 TP$ and $\psi: P\rightarrow\R^8$ a real analytic
embedding (immersion).  There exist $\epsilon >0$ and a unique
family $\{\psi_t:t\in (-\epsilon,\epsilon)\}$ of real analytic maps
$\psi_t:P\rightarrow\R^8$ with $\psi_0=\psi$ satisfying
\begin{equation}\label{ch5s1eq2}
\left(\frac{d\psi_t}{dt}\right)^e =
(\psi_t)_*(\chi)^{abc}(\Phi_0)_{abcd} (g_0)^{de}
\end{equation}
\noindent using index notation for tensors on $\R^8$, where
$(g_0)^{de}$ is the inverse of the Euclidean metric on $\R^8$.
Define $\Psi:(-\epsilon,\epsilon) \times P\rightarrow \R^8$ by
$\Psi(t,p)=\psi_t(p)$.  Then $M = \text{\emph{Image}}\, \Psi$ is a
nonsingular embedded (immersed) Cayley 4-fold in $\R^8$.
\end{thm}

\subsection{The symmetries method}\label{symmethod}

Now that we have a means of constructing calibrated submanifolds of
$\R^7$ and $\R^8$, we shall consider the situation where the
submanifold has a large symmetry group. The imposition of symmetry
on the system reduces its complexity.  This observation motivates
our method of construction, which is a generalisation of the work of
Joyce in \cite{Joy2}.

We know from the remarks after Definitions \ref{ch2s3dfn1} and
\ref{ch2s4dfn1} that it is natural to consider subgroups of $\G2
\ltimes \R^7$ or $\Spin(7)\ltimes\R^8$ as symmetry groups for our
calibrated submanifolds.

\medskip

Let us consider, for example, the associative case.  Suppose that
$\text{G}$ is a Lie subgroup of $\G2 \ltimes \R^7$ which has a
two-dimensional orbit $\mathcal{O}\subseteq\R^7$. Theorem
\ref{ch4s1thm3} allows us to evolve each point in $\mathcal{O}$
transversely to the action of $\text{G}$ and hence, hopefully,
construct an associative 3-fold with symmetry group $\text{G}$.

Formally, take $\chi$ to be a nowhere vanishing section of
$\Lambda^2T\text{G}$, which can easily be determined by finding a
basis for the Lie algebra of $\text{G}$.  Define
$\psi:\text{G}\rightarrow\mathcal{O}\subseteq\R^7$ to be an
embedding given by
$$\psi(\gamma)=\gamma\cdot(x_1,\ldots,x_7)$$
for $\gamma\in\text{G}$, where $(x_1,\ldots,x_7)$ is a point in
$\mathcal{O}$ and $\gamma\cdot(x_1,\ldots,x_7)$ denotes the action
of $\text{G}$ on $\R^7$.  Finally, for $t\in\R$, let
$\psi_t:\text{G}\rightarrow\R^7$ be given by
$$\psi_t(\gamma)=\gamma\cdot\big(x_1(t),\ldots,x_7(t)\big),$$
where $x_1(t),\ldots,x_7(t)$ are smooth real-valued functions of $t$
with $x_j(0)=x_j$ for $j=1,\ldots,7$.

We may thus calculate either side of \eq{ch4s1eq1} and get a coupled
system of seven first-order differential equations in seven
variables dependent on $t$; that is, of the form
$$\frac{d}{dt}\,\big(x_1(t),\ldots,x_7(t)\big)=\Big(y_1\big(x_1(t),\ldots,x_7(t)\big),\ldots,y_7\big(x_1(t),\ldots,
x_7(t)\big)\Big)$$ for functions $y_1,\ldots,y_7:\R^7\rightarrow\R$.

\begin{remark}
$y_1,\ldots,y_7$ are quadratic functions of their arguments.
\end{remark}

 By Theorem
\ref{ch4s1thm3}, a unique solution to this system exists for
$t\in(-\epsilon,\epsilon)$, for some $\epsilon>0$. Moreover, if
$$M=\big\{\gamma\cdot\big(x_1(t),\ldots,x_7(t)\big)\,:\,\gamma\in\text{G},\,t\in(-\epsilon
,\epsilon)\big\},$$ it is an associative 3-fold in $\R^7$ which is
clearly $\text{G}$-invariant.

\medskip

For the coassociative case, we need to consider Lie subgroups
$\text{G}$ of $\G2\ltimes\R^7$ which have a 3-dimensional orbit
$\mathcal{O}$.  However, we also need to choose $\mathcal{O}$ so
that $\varphi_0|_{\mathcal{O}}=0$; i.e. we need
$\psi:\text{G}\rightarrow\mathcal{O}$ to be an embedding such that
$\psi^*(\varphi_0)\equiv 0$ on $\text{G}$.

To construct Cayley examples with symmetries, we need to focus on
Lie subgroups of $\Spin(7)\ltimes\R^8$ that have 3-dimensional
orbits.

\begin{remark}
If we write the system of differential equations defining
coassociative or Cayley 4-folds with symmetries in the form
$$\frac{d\mathbf{x}}{dt}=\mathbf{y}(\mathbf{x}),$$
the components of $\mathbf{y}$ will be cubic functions of the
variables in $\mathbf{x}$.
\end{remark}

\medskip

The author has looked at a variety of different subgroups and has
derived systems of differential equations defining associative,
coassociative and Cayley submanifolds.  However, in the majority of
situations the author has been unsuccessful in solving the system.
In Sections \ref{s4} and \ref{s5} we present two important cases
which we have been able to solve.  Some additional scenarios where
the author has had less fortune are discussed in Section \ref{s6}.

\section{$\U(1)$-invariant associative cones}\label{s4}

In this section, we consider associative 3-folds which are invariant
both under an action of $\U(1)$ on the $\C^3$ component of
$\R^7\cong\R \oplus \C^3$ and under dilations.

\begin{dfn}\label{ch4s2dfn2}
Let $\R^+$ denote the group of positive real numbers under
multiplication.  The group action of $\R^+ \times \U(1)$ on
$\R^7\cong\R\oplus\C^3$ is given by, for some fixed $\alpha_1,
\alpha_2, \alpha_3 \in \R$,
\begin{equation*}
(x_1,z_1,z_2,z_3) \longmapsto (rx_1, \, re^{is\alpha_1}z_1,
 \, re^{is\alpha_2}z_2,  \, re^{is\alpha_3}z_3)
\qquad r>0,\, s\in\R. 
\end{equation*}
To ensure we have a $\U(1)$ action in $\G2$, we choose $\alpha_1$,
$\alpha_2$, $\alpha_3$ to be coprime integers satisfying $\alpha_1 +
\alpha_2 + \alpha_3 = 0$.
\end{dfn}

Define smooth maps $\psi_t:\R^+ \times \U(1) \rightarrow \R^7$ by
\begin{equation}\label{ch4s2subs2eq2}
\psi_t(r,e^{is}) = \big(rx_1(t), \, re^{is\alpha_1}z_1(t),
 \, re^{is\alpha_2}z_2(t),  \, re^{is\alpha_3}z_3(t)\big),
\end{equation}
where $x_1(t)$, $z_1(t)=x_2(t)+ix_3(t)$, $z_2(t)=x_4(t)+ix_5(t)$ and
$z_3(t)=x_6(t)+ix_7(t)$ are smooth functions of $t$.

Using \eq{ch4s2subs2eq2} we calculate the tangent vectors to the
group action given in Definition \ref{ch4s2dfn2}:
\begin{align*}
\mathbf{u} & =(\psi_t)_* \left( \frac{\partial}{\partial r} \right)
= \sum_{j=1}^{7} x_j \frac{\partial}{\partial x_j}\,\,\text{and}
\\[4pt]
\mathbf{v} & =  (\psi_t)_* \left( \frac{\partial}{\partial s}
\right)\\[4pt] &=
  \alpha_1 \left( x_2 \frac{\partial}{\partial x_3} - x_3 \frac{\partial}{\partial x_2} \right)
  + \alpha_2 \left( x_4 \frac{\partial}{\partial x_5} - x_5 \frac{\partial}{\partial x_4}\right)
  + \alpha_3 \left( x_6 \frac{\partial}{\partial x_7} - x_7 \frac{\partial}{\partial x_6} \right).
\end{align*}

\noindent If we take $\chi = \frac{\partial}{\partial r} \w
\frac{\partial}{\partial s}$, then $(\psi_t)_* (\chi) = \mathbf{u} \w
\mathbf{v}$.  We deduce that, writing $\mathbf{e}_j=\frac
{\partial}{\partial x_j}$,
\begin{align*}
\mathbf{u}^a \mathbf{v}^b (\varphi_0)_{abc}(g_0)^{cd}
=&\,\,\big(\alpha_1(x_2^2 + x_3^2) + \alpha_2(x_4^2 + x_5^2)
+\alpha_3(x_6^2 + x_7^2)\big)\mathbf{e}_1 \\
&+ \big(-\alpha_1 x_1x_2 + (\alpha_2 -\alpha_3) (x_4x_7 +
x_5x_6)\big)\mathbf{e}_2\\  &+ \big(-\alpha_1 x_1x_3 + (\alpha_2
-\alpha_3) (x_4x_6 - x_5x_7)\big)\mathbf{e}_3\\& + \big(-\alpha_2
x_1x_4 + (\alpha_3 -\alpha_1) (x_2x_7+x_3x_6)\big)\mathbf{e}_4\\& +
\big(-\alpha_2 x_1x_5 + (\alpha_3 -\alpha_1) (x_2x_6 -
x_3x_7)\big)\mathbf{e}_5\\&+ \big(-\alpha_3 x_1x_6 + (\alpha_1
-\alpha_2) (x_2x_5 + x_3x_4)\big)\mathbf{e}_6\\& + \big(-\alpha_3
x_1x_7 + (\alpha_1 -\alpha_2) (x_2x_4 - x_3x_5)\big)\mathbf{e}_7.
\end{align*}
\noindent We also have that
\begin{equation*}
\frac{d\psi_t}{dt}  =
\sum_{j=1}^{7}\frac{dx_j(t)}{dt}\,\mathbf{e}_j.
\end{equation*}

Equating both sides of \eq{ch4s1eq1} using the above formulae as
described in $\S$\ref{symmethod}, we obtain the following theorem.

\begin{thm}\label{ch4s2thm2} Use the notation of Definition
\ref{ch4s2dfn2}.  Let $\beta_1 = \alpha_2 - \alpha_3$, $\beta_2 =
\alpha_3 - \alpha_1$ and $\beta_3 = \alpha_1 - \alpha_2$.  Let
$x_1(t)$ be a smooth real-valued function of $t$ and let $z_1(t)$,
$z_2(t)$, $z_3(t)$ be smooth complex-valued functions of $t$ such
that
\begin{align}
\label{ch4s2subs2eq3}
\frac{dx_1}{dt} & =  \alpha_1|z_1|^2 + \alpha_2|z_2|^2 +
\alpha_3|z_3|^2,
\\[4pt]
\label{ch4s2subs2eq4}
\frac{dz_1}{dt} & =  -\alpha_1x_1z_1 + i \beta_1 \overline{z_2z_3},
\\[4pt]
\label{ch4s2subs2eq5}
\frac{dz_2}{dt} & =  -\alpha_2x_1z_2 + i \beta_2
\overline{z_3z_1}\,\,\text{and}
\\[4pt]
\label{ch4s2subs2eq6}
\frac{dz_3}{dt} & =  -\alpha_3x_1z_3 + i \beta_3 \overline{z_1z_2}.
\end{align}
These equations have a solution for all $t\in\R$ and the subset $M$
of\/ $\R\oplus\C^3\cong\R^7$ defined by
\begin{equation*}
M  =
\big\{\big(rx_1(t),\,re^{is\alpha_1}z_1(t),\,re^{is\alpha_2}z_2(t) ,
\,re^{is\alpha_3}z_3(t)\big):\,r\in\R^+,\,s,t\in\R\big\}
\end{equation*}
is an associative 3-fold in $\R^7$.  Moreover,
\eq{ch4s2subs2eq3}-\eq{ch4s2subs2eq6} imply that
$x_1^2+|z_1|^2+|z_2|^2+|z_3|^2$ can be chosen to be $1$ and that
$\Re(z_1z_2z_3)=A$, where $A$ is a real constant.
\end{thm}

\begin{proof} Noting that $\beta_1 + \beta_2 + \beta_3 = 0$, we immediately
see that $x_1^2 + |z_1|^2 + |z_2|^2 + |z_3|^3$ is a constant which
we can take to be one.  This is hardly surprising since the
associative 3-fold was constructed so as to be a cone.  We also see
from \eq{ch4s2subs2eq4}-\eq{ch4s2subs2eq6} that
\begin{equation*}
\frac{d}{dt}\,(z_1z_2z_3)  =
i(\beta_1|z_2|^2|z_3|^2+\beta_2|z_3|^2|z_1|^2+\beta_3|z_1|^2|z_2|^2),
\end{equation*}
which is purely imaginary, and therefore $\Re(z_1z_2z_3)=A$ is a
 constant.

Notice that the functions $x_1$, $z_1$, $z_2$ and $z_3$ are bounded,
hence their first derivatives are bounded by
\eq{ch4s2subs2eq3}-\eq{ch4s2subs2eq6}.  Thus, all of the functions
which determine the behaviour of the solutions to
\eq{ch4s2subs2eq3}-\eq{ch4s2subs2eq6} are bounded, from which it
follows that they have solutions for all $t\in\R$.
\end{proof}

Writing $z_j(t)=r_j(t)e^{i\theta_j(t)}$ for $j=1,2,3$ and $\theta =
\theta_1 + \theta_2 + \theta_3$,
\eq{ch4s2subs2eq3}-\eq{ch4s2subs2eq6} become
\begin{align}
\label{ch4s2subs2eq7}
\frac{dx_1}{dt} & =  \alpha_1 r_1^2 + \alpha_2 r_2^2 + \alpha_3
r_3^2;
\\[4pt]
\label{ch4s2subs2eq8}
\frac{dr_1}{dt} & =  -\alpha_1 x_1 r_1 + \beta_1 r_2 r_3 \sin
\theta;
\\[4pt]
\label{ch4s2subs2eq9}
\frac{dr_2}{dt} & =  -\alpha_2 x_1 r_2 + \beta_2 r_3 r_1 \sin
\theta;
\\[4pt]
\label{ch4s2subs2eq10}
\frac{dr_3}{dt} & =  -\alpha_3 x_1 r_3 + \beta_3 r_1 r_2 \sin
\theta;\,\text{and}
\\[4pt]
\label{ch4s2subs2eq11}
 r_j^2 \frac{d\theta_j}{dt} & =  \beta_j A
\qquad\qquad\qquad\qquad \quad\!\!\text{for $j=1,2,3$,}
\end{align}
with the conditions
\begin{align}
\label{ch4s2subs2eq12}
x_1^2 + r_1^2 + r_2^2 + r_3^2 & =
1\,\;\text{and}\\
\label{ch4s2subs2eq13}
r_1r_2r_3\cos\theta & =  A.
\end{align}

\medskip

We notice that we are restricted in our choices of the real
parameter $A$. The problem of maximising $A^2$, by
\eq{ch4s2subs2eq12} and \eq{ch4s2subs2eq13}, is equivalent to the
problem of maximising $r_1^2r_2^2r_3^2$ subject to $r_1^2 + r_2^2 +
r_3^2 = 1$.  By direct calculation the solution\\[4pt] is $r_1^2 = r_2^2
= r_3^2 = \frac{1}{3}$. Therefore $A \in \left[-\frac{1}{3
\sqrt{3}}, \frac{1}{3 \sqrt{3}} \right]$.  We can restrict to
$A\geq0$\linebreak\\[-12pt] since changing the sign of $A$ corresponds
to reversing the sign of $\cos \theta$, so the addition of $\pi$ to
$\theta$.

The case where $A = \frac{1}{3 \sqrt{3}}$ is immediately soluble
since this forces $r_1 = r_2 = r_3 = \frac{1}{\sqrt{3}}\,$, which
implies $x_1 = 0$ by \eq{ch4s2subs2eq12} and $\cos \theta = 1$ by
\eq{ch4s2subs2eq13}, so we can take $\theta = 0$. Equations
\eq{ch4s2subs2eq11} become
\begin{equation*}
\begin{split}
\frac{1}{3} \frac{d\theta_j}{dt} = \frac{1}{3 \sqrt{3}}\,\beta_j
\qquad \text{for $j=1,2,3$,}
\end{split}
\end{equation*}
\noindent which can easily be solved, along with the condition
$\theta = 0$, to give:
\begin{equation*}
\theta_j(t) = \frac{\beta_j}{\sqrt{3}}\, t + \gamma_j
\qquad\text{for $j=1,2,3,$}
\end{equation*}
\noindent where $\gamma_1$, $\gamma_2$, $\gamma_3$ are real
constants which sum to zero.  Then
\begin{equation*}
M = \left\{\left(0, \hspace{2pt} r e^{i\phi_1}, \hspace{2pt} r
e^{i\phi_2}, \hspace{2pt} r e^{i\phi_3} \right): r>0, \hspace{2pt}
\phi_1, \phi_2, \phi_3 \in \R, \hspace{2pt} \phi_1 + \phi_2 + \phi_3
= 0  \right\},
\end{equation*}
which is a $\U(1)^2$-invariant special Lagrangian cone, as studied
in \cite[$\S$III.3.A]{HarLaw}, embedded in $\R^7$ and is therefore
in itself not an interesting object of study here.  Any associative
3-fold constructed with $x_1 = 0$ will be at least a
$\U(1)$-invariant special Lagrangian cone and so we shall not
consider this situation further.  However, we know that $M$ must be
the limiting case of the family of associative 3-folds parameterised
by $A$ as it tends to $\frac{1}{3 \sqrt{3}}\,$.

\medskip

We may also solve the equations in the following special case.

\begin{thm}\label{ch4s2thm3}
Use the notation of Theorem \ref{ch4s2thm2}.  Suppose that
$\alpha_2=\alpha_3$.  Then $x_1,z_1,z_2$ and $z_3$ may be chosen to
satisfy $x_1^2+|z_1|^2+|z_2|^2+|z_3|^2=1$ and\/ $\Im z_1=0$.
Moreover, they satisfy:
\begin{gather*}\Re(z_1z_2z_3)=A;\qquad
|z_1|(x_1^2+|z_1|^2-1)=B;
\\
\Re(z_1(z_2^2-z_3^2))=C;\quad\text{and}\quad
\Im(z_1(z_2^2+z_3^2))=D
\end{gather*}
for some real constants $A$, $B$, $C$ and $D$.
\end{thm}

\begin{proof}
Since $\beta_1=0$, \eq{ch4s2subs2eq11} implies that the argument of
$z_1$ is constant.  Using $\U(1)$ we can take it to be zero so that
$z_1$ is real.  Moreover, $\beta_1=0$ and \eq{ch4s2subs2eq12} imply
that $x_1$ and $z_1$ evolve amongst themselves and hence, using
\eq{ch4s2subs2eq3} and \eq{ch4s2subs2eq4}, we deduce that the real
function $f=|z_1|(x_1^2+|z_1|^2-1)$ is constant.  Note that $\SU(2)$
acts on the $(z_2,z_3)$-plane.  We are thus led to calculate
\begin{align*}
\frac{d}{dt}\,\big(z_1(az_2+bz_3)&(-\bar{b}z_2+\bar{a}z_3)\big)\\
&=-4i\beta_2|z_1|^2\Re(a\bar{b}z_2z_3)+i\beta_2|z_1|^2(|a|^2-|b|^2)(|z_3|^2-|z_2|^2)\end{align*}
for constants $a,b\in\C$, which is purely imaginary.  Equating real
parts for $(a,b)=(1,-1)$ and $(a,b)=(i,1)$ leads to the final two
conserved quantities in the statement of the theorem.
\end{proof}

In Theorem \ref{ch4s2thm3}, we have six conditions on seven variables,
which thus determine the solution to the system of differential
equations \eq{ch4s2subs2eq3}-\eq{ch4s2subs2eq6} and hence the
associative cone constructed by Theorem \ref{ch4s2thm2} for
$\alpha_2=\alpha_3$. Moreover, we may construct a function
$\pi:\R\oplus\C^3\rightarrow\R^6$ by mapping $(x_1,z_1,z_2,z_3)$ to
the six real constant functions given in Theorem \ref{ch4s2thm3},
which are defined by the initial values
$\big(x_1(0),z_1(0),z_2(0),z_3(0)\big)$.

\emph{Sard's Theorem} \cite[p. 173]{Lang} states that if
$f:X\rightarrow Y$ is a smooth map between finite-dimensional
manifolds, the set of $y\in Y$ with some $x\in f^{-1}(y)$ such that
$df|_x:T_xX\rightarrow T_yY$ is \emph{not} surjective is of measure
zero in $Y$. Therefore, $f^{-1}(y)$ is a submanifold of $X$ of
dimension $\text{dim}\,X-\text{dim}\,Y$ for almost all $y\in Y$.
Applying Sard's Theorem, generically the fibres of $\pi$ will be
1-dimensional submanifolds of $\R\oplus\C^3\cong\R^7$. Moreover, we
know that these fibres are compact by the conditions in Theorem
\ref{ch4s2thm3}. Hence, the variables form loops in $\R^7$ for
generic initial values; i.e. the solutions are \emph{periodic} in
$t$.  We deduce the following result.

\begin{thm}\label{ch4s2thm4}  Use the notation of Theorem
\ref{ch4s2thm2} and suppose that $\alpha_2=\alpha_3$.  For generic
values of the functions $x_1$, $z_1$, $z_2$ and $z_3$ at $t=0$, the
associative 3-folds constructed by Theorem \ref{ch4s2thm2} are
closed $\U(1)$-invariant cones over $T^2$ in $\R^7$.
\end{thm}

This family of cones is determined by four real parameters, whereas
the corresponding SL family, as discussed in \cite[$\S$7]{Joy2}, is
parameterised by one rational variable.  Therefore, these cones are
generically not SL.

\medskip

We may also apply the theory described in \cite[$\S$6]{Lotay1} to
the family of cones given in Theorem \ref{ch4s2thm4} to produce
examples of \emph{ruled} associative 3-folds which are
\emph{asymptotically conical}.  We thus define the terms we require,
noting that a cone $C$ in $\R^7$ is said to be \emph{two-sided} if
$C=-C$.

\begin{dfn}
\label{ruled} Let $M$ be a 3-dimensional submanifold of $\R^7$.  A
\emph{ruling} of $M$ is a pair $(\Sigma,\pi)$, where $\Sigma$ is a
2-dimensional manifold and $\pi:M \rightarrow\Sigma$ is a smooth
map, such that for all $\sigma\in\Sigma$ there exist ${\bf
v}_{\sigma}\in \mathcal{S}^6$, ${\bf w}_{\sigma}\in\R^7$ such that
$\pi^{-1}(\sigma) =\{r{\bf v}_{\sigma}+{\bf w}_{\sigma}:r\in\R\}$.
Then the triple $(M,\Sigma, \pi)$ is a \emph{ruled submanifold} of
$\R^7$.

An \emph{r-orientation} for a ruling $(\Sigma,\pi)$ of $M$ is a
choice of orientation for the affine straight line
$\pi^{-1}(\sigma)$ in $\R^7$, for each $\sigma\in\Sigma$, which
varies smoothly with $\sigma$.  A ruled submanifold with an
r-orientation for the ruling is called an \emph{r-oriented ruled
submanifold}.

Let $(M,\Sigma,\pi)$ be an r-oriented ruled submanifold. For each
$\sigma \in\Sigma$, let $\phi(\sigma)$ be the unique unit vector in
$\R^7$ parallel to $\pi^{-1}(\sigma)$ and in the positive direction
with respect to the orientation on $\pi^{-1}(\sigma)$, given by the
r-orientation. Then $\phi:\Sigma\rightarrow \mathcal{S}^6$ is a
smooth map. Define $\psi:\Sigma \rightarrow\R^7$ such that, for all
$\sigma\in\Sigma$, $\psi(\sigma)$ is the unique vector in
$\pi^{-1}(\sigma)$ orthogonal to $\phi(\sigma)$.  Then $\psi$ is a
smooth map and we may write
\begin{equation}
\begin{split}
\label{ruled1}
M=\{r\phi(\sigma)+\psi(\sigma):\sigma\in\Sigma,\hspace{2pt}
r\in\R\}.
\end{split}
\end{equation}

Define the \emph{asymptotic cone} $M_0$ of a ruled submanifold $M$
by
\begin{equation*}
\begin{split}
M_0=\{{\bf v}\in\R^7:\text{${\bf v}$ is parallel to
$\pi^{-1}(\sigma)$ for some $\sigma\in\Sigma$}\}.
\end{split}
\end{equation*}
If $M$ is also r-oriented, then
\begin{equation}
\begin{split}
\label{ruled2}
M_0=\{r\phi(\sigma):\sigma\in\Sigma,\hspace{2pt}r\in\R\}
\end{split}
\end{equation}
\noindent and is usually a 3-dimensional two-sided cone; that is,
whenever $\phi$ is an immersion.
\end{dfn}

\begin{dfn}
\label{asym} Let $M_0$ be a closed cone in $\R^7$ and let $M$ be a
closed nonsingular submanifold in $\R^7$.  We say that $M$ is
\emph{asymptotically conical to $M_0$ with rate $\alpha$}, for some
$\alpha<1$, if there exists some constant $R>0$, a compact subset
$K$ of $M$ and a diffeomorphism
$\Psi:M_0\setminus\bar{B}_R\rightarrow M\setminus K$ such that
\begin{equation*}
\big|\nabla^{k}\big(\Psi({\bf x})-\iota({\bf x})\big)\big|=O\big(r^{\alpha-k}\big) \quad
\text{for $k\in\N$ as $r\rightarrow\infty$,}
\end{equation*}
\noindent where $\bar{B}_R$ is the closed ball of radius $R$ in
$\R^7$ and $\iota
:M_0\rightarrow\R^7$ is the inclusion map.  Here
$|\,.\,|$ is calculated using the cone metric on
$M_0\setminus\bar{B}_R$, and $\nabla$ is a combination of the
Levi--Civita connection derived from the cone metric and the flat
connection on $\R^n$, which acts as partial differentiation.
\end{dfn}

We now use the construction involving \emph{holomorphic vector
fields} given in \cite[Proposition 6.8]{Lotay1}

\begin{thm} Use the notation of Theorem \ref{ch4s2thm2} and suppose
that $\alpha_2=\alpha_3=-1$.  Let $M$, as given in Theorem
\ref{ch4s2thm2}, be an associative cone over $T^2$, which, by
Theorem \ref{ch4s2thm4}, occurs for generic choices of $x_1(0)$,
$z_1(0)$, $z_2(0)$ and $z_3(0)$.  Let $u,v:\R^2\rightarrow\R$ be
functions satisfying the Cauchy--Riemann equations and let
$M_0=M\cup (-M)\cup\{0\}$. The subset $M_{u,v}$ of\/ $\R\oplus\C^3$
given by
\begin{align*}M_{u,v}=\Big\{\Big(&rx_1(t)+v(s,t)\big(2|z_1(t)|^2-|z_2(t)|^2-|z_3(t)|^2\big),\\
&\,e^{2is}\big(r+2iu(s,t)-2v(s,t)x_1(t)\big)z_1(t),\\&\,e^{-is}\big(\big(r-iu(s,t)+v(s,t)x_1(t)
\big)z_2(t)
-3iv(s,t)\overline{z_3z_1}\,\big),\\
&\,e^{-is}\big(\big(r-iu(s,t)+v(s,t)x_1(t)\big)z_3(t)+3iv(s,t)\overline{z_1z_2}\,\big)
\Big):\,r,s,t\in\R\Big\}\end{align*} is an r-oriented ruled
associative 3-fold in $\R^7\cong\R\oplus\C^3$.  Moreover, $M_{u,v}$
is asymptotically conical to $M_0$ with rate $-1$ in the sense of
Definition \ref{asym}.
\end{thm}

\begin{proof}
Define $\phi:\R^2\rightarrow\R\oplus\C^3\cong\R^7$ by
$$\phi(s,t)=(x_1(t),\,e^{2is}z_1(t),\,e^{-is}z_2(t),\,e^{-is}z_3(t)).$$
Since $x_1^2+|z_1|^2+|z_2|^2+|z_3|^2=1$, $\phi$ maps into
$\mathcal{S}^6$, and we can write $M_0$ in the form \eq{ruled2}.
Define a holomorphic vector field $w$ using $u$ and $v$ as follows:
\begin{equation*}
w = u(s,t)\frac{\partial}{\partial s} +
v(s,t)\frac{\partial}{\partial t}\,.
\end{equation*}
Define $\psi=\mathcal{L}_w\phi$, where $\mathcal{L}_w$ denotes the
Lie deriviative with respect to $w$, and define $M_{u,v}$ by
\eq{ruled1} for these choices of $\phi$ and $\psi$. Calculation
using equations \eq{ch4s2subs2eq3}-\eq{ch4s2subs2eq6} of Theorem
\ref{ch4s2thm2} shows that $M_{u,v}$ can be written as stated in the
theorem. Applying \cite[Proposition 6.8 \& Theorem 6.9]{Lotay1},
since $M_0$ is a cone over $T^2$, gives us the various properties of
$M_{u,v}$ as claimed.
\end{proof}

\begin{remark} Although $M$ and hence $M_0$ is $\U(1)$-invariant, $M_{u,v}$
will not be in general. \end{remark}

\section{$\SU(2)$-invariant Cayley 4-folds}\label{s5}

We consider three different natural actions of $\SU(2)$ on
$\C^4\cong\R^8$ in $\Spin(7)$, though the first two only give
trivial examples of Cayley 4-folds. The first is where $\SU(2)$ acts
on $\C^4\cong\C^2\oplus\C^2$ in the usual manner upon one $\C^2$ and
trivially upon the other. The construction using this action gives
an affine $\C^2\subseteq\C^4$ as the Cayley 4-fold.  The second is
where $\SU(2)$ acts on $\C^4\cong\C^3\oplus\C$ as $\SO(3)$ on $\C^3$
and trivially on $\C$. The construction then produces a complex
surface in $\C^4$ as the Cayley 4-fold, which may be written as
follows:
\begin{equation*}
\left\{(z_1,z_2,z_3,z_4):z_1^2+z_2^2+z_3^2=A, \,
z_4=B\right\},\,\;\text{where $A,B\in\C$ are constants}.
\end{equation*}
We therefore turn our attention to the diagonal action of $\SU(2)$.

\begin{dfn}\label{ch5s3subs2dfn1} Let
\begin{equation*}
X = \left(\begin{array}{rr} a & b \\ -\bar{b} &
\bar{a}\end{array}\right)\in\SU(2),
\end{equation*}
where $a,b\in\C$ such that $|a|^2+|b|^2=1$.  Then $X$ acts on
$(z_1,z_2,z_3,z_4)\in\C^4\cong\R^8$ as:
\begin{equation*}
X\cdot(z_1,z_2,z_3,z_4) =
(az_1+bz_2,-\bar{b}z_1+\bar{a}z_2,az_3+bz_4,
-\bar{b}z_3+\bar{a}z_4).
\end{equation*}
\end{dfn}

Define smooth maps $\psi_t:\SU(2)\rightarrow\C^4\cong\R^8$ by:
\begin{equation*}
\psi_t(X) = X\cdot\big(z_1(t),z_2(t),z_3(t),z_4(t)\big),
\end{equation*}
where $z_1(t)$, $z_2(t)$, $z_3(t)$ and $z_4(t)$ are smooth functions
of $t$.

Calculation shows that we may take the following three complex
matrices as a basis for the Lie algebra of $\SU(2)$ acting in this
way:
\begin{gather*}
U_1 = \left(\begin{array}{cccc} i & 0 & 0 & 0 \\ 0 & -i & 0 & 0 \\
0 & 0 & i & 0 \\ 0 & 0 & 0 & -i \end{array}\right); \qquad
U_2 = \left(\begin{array}{cccc} 0 & 1 & 0 & 0 \\ -1 & 0 & 0 & 0 \\
0 & 0 & 0 & 1 \\ 0 & 0 & -1 & 0 \end{array}\right);\,\\[4pt]
\quad\text{and}\quad
U_3 = \left(\begin{array}{cccc} 0 & i & 0 & 0 \\ i & 0 & 0 & 0 \\
0 & 0 & 0 & i \\ 0 & 0 & i & 0 \end{array}\right).
\end{gather*}
\noindent If we let $\mathbf{u}_j=(\psi_t)_*(U_j)$ for $j=1,2,3$,
\begin{align*}
\mathbf{u}_1 & =  i\!\left(z_1\frac{\partial}{\partial
z_1}-\bar{z}_1
\frac{\partial}{\partial\bar{z}_1}-z_2\frac{\partial}{\partial z_2}+
\bar{z}_2\frac{\partial}{\partial\bar{z}_2}+z_3
\frac{\partial}{\partial z_3}-\bar{z}_3
\frac{\partial}{\partial\bar{z}_3}-z_4\frac{\partial}{\partial z_4}+
\bar{z}_4\frac{\partial}{\partial\bar{z}_4}\right), \\[4pt]
\mathbf{u}_2 & =  z_2\frac{\partial}{\partial z_1}+\bar{z}_2
\frac{\partial}{\partial\bar{z}_1}-z_1\frac{\partial}{\partial z_2}-
\bar{z}_1\frac{\partial}{\partial\bar{z}_2}+z_4
\frac{\partial}{\partial z_3}+\bar{z}_4
\frac{\partial}{\partial\bar{z}_3}-z_3\frac{\partial}{\partial z_4}-
\bar{z}_3\frac{\partial}{\partial\bar{z}_4}\,\;\,\text{and} \\[4pt]
\mathbf{u}_3 & = i\!\left(z_2\frac{\partial}{\partial z_1}-\bar{z}_2
\frac{\partial}{\partial\bar{z}_1}+z_1\frac{\partial}{\partial z_2}-
\bar{z}_1\frac{\partial}{\partial\bar{z}_2}+z_4
\frac{\partial}{\partial z_3}-\bar{z}_4
\frac{\partial}{\partial\bar{z}_3}+z_3\frac{\partial}{\partial z_4}-
\bar{z}_3\frac{\partial}{\partial\bar{z}_4}\right).
\end{align*}
Thus, if we take $\chi=U_1\w U_2\w U_3$, $(\psi_t)_*(\chi)
=\mathbf{u}_1\w \mathbf{u}_2\w \mathbf{u}_3$.  Using the equations
above for $\mathbf{u}_j$ and the formula \eq{ch2s4eq1} for $\Phi_0$,
we may calculate the right-hand side of \eq{ch5s1eq2}:
\begin{align*}
\mathbf{u}_1^a\mathbf{u}_2^b\mathbf{u}_3^c&(\Phi_0)_{abcd}(g_0)^{de}
\\
&=\left(z_1\left(|z_1|^2+|z_2|^2+|z_3|^2-|z_4|^2\right)+2(\overline{z_1z_4-
z_2z_3}+z_2z_3)\bar{z}_4\right)\frac{\partial}{\partial z_1}\\
&{}+\left(\bar{z}_1\left(|z_1|^2+|z_2|^2+|z_3|^2-|z_4|^2\right)+
2(z_1z_4-z_2z_3+\overline{z_2z_3})z_4\right)\frac{\partial}{\partial
\bar{z}_1} \\
&{}+\left(z_2\left(|z_1|^2+|z_2|^2-|z_3|^2+|z_4|^2\right)-2(\overline{z_1z_4-
z_2z_3}-z_1z_4)\bar{z}_3\right)\frac{\partial}{\partial z_2} \\
&{}+\left(\bar{z}_2\left(|z_1|^2+|z_2|^2-|z_3|^2+|z_4|^2\right)-
2(z_1z_4-z_2z_3-\overline{z_1z_4})z_3\right)\frac{\partial}{\partial
\bar{z}_2} \displaybreak[0]\\
&{}+\left(z_3\left(|z_1|^2-|z_2|^2+|z_3|^2+|z_4|^2\right)-2(\overline{z_1z_4-
z_2z_3}-z_1z_4)\bar{z}_2\right)\frac{\partial}{\partial z_3} \\
&{}+\left(\bar{z}_3\left(|z_1|^2-|z_2|^2+|z_3|^2+|z_4|^2\right)-
2(z_1z_4-z_2z_3-\overline{z_1z_4})z_2\right)\frac{\partial}{\partial
\bar{z}_3} \\
&{}+\left(z_4\left(-|z_1|^2+|z_2|^2+|z_3|^2+|z_4|^2\right)+2(\overline{z_1z_4-
z_2z_3}+z_2z_3)\bar{z}_1\right)\frac{\partial}{\partial z_4}\\
&{}+ \left(\bar{z}_4\left(-|z_1|^2+|z_2|^2+|z_3|^2+|z_4|^2\right)+
2(z_1z_4-z_2z_3+\overline{z_2z_3})z_1\right)\frac{\partial}{\partial
\bar{z}_4}\,.
\end{align*}
Moreover,
\begin{equation*}
\frac{d\psi_t}{dt} = \sum_{j=1}^4
\frac{dz_j}{dt}\frac{\partial}{\partial z_j} + \sum_{j=1}^4
\frac{d\bar{z}_j}{dt} \frac{\partial}{\partial\bar{z}_j}.
\end{equation*}
\noindent Equating both sides of \eq{ch5s1eq2} and using Theorem
\ref{ch5s1thm5} gives the following result.

\begin{thm}\label{ch5s3subs2thm1}
Let $z_1(t),z_2(t),z_3(t),z_4(t)$ be smooth complex-valued functions
of\/ $t$ satisfying
\begin{align}\label{ch5s3subs2eq4}
\frac{dz_1}{dt} & = z_1\left(|z_1|^2+|z_2|^2+|z_3|^2-|z_4|^2\right)+
2(\overline{z_1z_4-z_2z_3}+z_2z_3)\bar{z}_4, \\[4pt]
\frac{dz_2}{dt} &= z_2\left(|z_4|^2+|z_1|^2+|z_2|^2-|z_3|^2\right)-
2(\overline{z_1z_4-z_2z_3}-z_1z_4)\bar{z}_3, \\[4pt]
\frac{dz_3}{dt} &= z_3\left(|z_3|^2+|z_4|^2+|z_1|^2-|z_2|^2\right)-
2(\overline{z_1z_4-z_2z_3}-z_1z_4)\bar{z}_2\,\,\text{and} \\[4pt]
\label{ch5s3subs2eq7}
\frac{dz_4}{dt} & = z_4\left(|z_2|^2+|z_3|^2+|z_4|^2-|z_1|^2\right)+
2(\overline{z_1z_4-z_2z_3}+z_2z_3)\bar{z}_1
\end{align}
for all $t\in (-\epsilon,\epsilon)$, for some $\epsilon>0$.  The
subset $M$ of\/ $\C^4\cong\R^8$ defined by
\begin{equation*}
M = \big\{X\cdot\big(z_1(t),z_2(t),z_3(t),z_4(t)\big)\,:\,t\in
(-\epsilon,\epsilon),\, X\in\SU(2)\big\},
\end{equation*}
\noindent where the action of\/ $\SU(2)$ on $\C^4$ is given in
Definition \ref{ch5s3subs2dfn1}, is a Cayley 4-fold in $\R^8$.
\end{thm}

We are able to give an explicit description of the Cayley 4-folds
constructed in Theorem \ref{ch5s3subs2thm1}. 
Let $u(t)$ be a real-valued function satisfying
\begin{equation}\label{ch5s3subs2eq8}
\frac{du}{dt} = 2(|z_1|^2+|z_2|^2+|z_3|^2+|z_4|^2)u.
\end{equation}
We observe, using \eq{ch5s3subs2eq4}-\eq{ch5s3subs2eq7}, that the
following quadratics satisfy \eq{ch5s3subs2eq8}:
\begin{align*}
&|z_1|^2-|z_2|^2+|z_3|^2-|z_4|^2;   &z_1\bar{z}_2+z_3\bar{z}_4; \\
&{\rm Re}(z_1z_4-z_2z_3);\quad\text{and} &z_1\bar{z}_3+z_2\bar{z}_4.
\end{align*}
\noindent Hence, each of these quadratics is a constant multiple of
$u$. The first two correspond to the moment maps of the $\SU(2)$
action and the latter two are $\SU(2)$-invariant.  The first two
quadratics are not $\SU(2)$-invariant, but
\begin{align}
 Q&(z_1,z_2,z_3,z_4) =
(|z_1|^2-|z_2|^2+|z_3|^2-|z_4|^2)^2+4|z_1\bar{z}_2+z_3
\bar{z}_4|^2\nonumber\\
  &=  (|z_1|^2+|z_2|^2)^2+(|z_3|^2+|z_4|^2)^2+
  2|z_1\bar{z}_3+z_2\bar{z}_4|^2-2|z_1z_4-z_2z_3|^2
\label{ch5s3subs2eq9}
\end{align}
\noindent is $\SU(2)$-invariant and is a constant multiple of $u^2$.

Using \eq{ch5s3subs2eq4}-\eq{ch5s3subs2eq7}, we calculate
\begin{equation*}
\frac{d}{dt}\Im(z_1z_4-z_2z_3)=-2
(|z_1|^2+|z_2|^2+|z_3|^2+|z_4|^2)\Im(z_1z_4-z_2z_3).
\end{equation*}
\noindent Therefore, by \eq{ch5s3subs2eq8}, $\Im(z_1z_4-z_2z_3)$ is
a constant multiple of $u^{-1}$ and is an $\SU(2)$-invariant
quadratic. We then state our result, which is immediate from
our discussion above.


\begin{thm}\label{ch5s3subs2thm2}
Let $A$, $B$, $C$ and $D$ be real constants.  Let
$M\subseteq\C^4\cong\R^8$ be defined by
\begin{equation*}
M = \{X\cdot(z_1,z_2,z_3,z_4):X\in\SU(2)\},
\end{equation*}
\noindent where the action of $X\in\SU(2)$ on $\C^4$ is given in
Definition \ref{ch5s3subs2dfn1} and $z_1,z_2,z_3,z_4$ satisfy:
\begin{align}
\label{ch5s3subs2eq12}
Q(z_1,z_2,z_3,z_4)\,\big(\Im(z_1z_4-z_2z_3)\big)^2 & =  A; \\
\Re(z_1z_4-z_2z_3)\,\Im(z_1z_4-z_2z_3) & =  B; \\
\Re(z_1\bar{z}_3+z_2\bar{z}_4)\,\Im
(z_1z_4-z_2z_3) & =  C;\,\text{and} \\
\label{ch5s3subs2eq15}
\Im(z_1\bar{z}_3+z_2\bar{z}_4)\,\Im (z_1z_4-z_2z_3) & =  D,
\end{align}
\noindent with $Q(z_1,z_2,z_3,z_4)$ given by \eq{ch5s3subs2eq9}.
Then $M$ is a Cayley 4-fold in $\R^8$.
\end{thm}
\noindent The set of conditions
\eq{ch5s3subs2eq12}-\eq{ch5s3subs2eq15} on the complex functions
$z_1,z_2,z_3,z_4$ consists of setting one real octic and three real
quartics to be constant, which defines a 4-dimensional subset of
$\C^4$.  Hence, Theorem \ref{ch5s3subs2thm2} completely describes
the $\SU(2)$-invariant Cayley 4-folds given by Theorem
\ref{ch5s3subs2thm1}.

\section{Further examples}\label{s6}

In this final section we present an example of a symmetry group and
its corresponding system of ordinary differential equations for each
type of calibrated submanifold considered in this paper. These
equations are derived using the method introduced in
$\S$\ref{symmethod}.  Since the calculations involved in this method
have already been described in detail through the work of the
previous two sections, we feel justified in our omission of the
relevant calculations here.

Though the author has had little success in attempting to solve the
systems in this section himself, it is hoped that their exposition
will be useful to others.

\subsection{Associative 3-folds invariant under a subgroup of\\
$\R\times\U(1)^2$}

We may decompose $\R^7 \cong \R \oplus \C^3$, and so the action of
$\R \times \U(1)^2$ on $\R^7$ may be written as:
\begin{equation}\label{ch4s2eq1}
(x_1, z_1, z_2, z_3) \longmapsto (x_1+c,\, e^{i\phi_1}z_1,\,
e^{i\phi_2}z_2,\, e^{-i(\phi_1 + \phi_2)}z_3), \quad \text{$c,
\phi_1, \phi_2 \in
\R$.} 
\end{equation}
However, we want a two-dimensional orbit, so we choose a
two-dimensional subgroup of $\R \times \U(1)^2$.

\begin{dfn}\label{ch4s2dfn1} Let $\lambda$, $\mu$, $\nu$ be real numbers which are not all
zero. Define $\text{G}$ to be the subgroup of $\R \times \U(1)^2$
which acts as in \eq{ch4s2eq1} with the following imposed:
\begin{equation}\label{ch4s2eq2}
\lambda c + \mu \phi_1 + \nu \phi_2 = 0. 
\end{equation}
If $\mu=\nu=0$, then $\text{G}$ is $\U(1)^2$.  Suppose $\mu\nu\neq
0$.  If there exist coprime integers $p$ and $q$ such that $\mu p +
\nu q = 0$, then $\text{G}$ is $\R\times\U(1)$ and otherwise it is an
$\R^2$ subgroup.
\end{dfn}

Using the method of $\S$\ref{symmethod} provides the following
theorem.

\begin{thm}\label{ch4s2thm1}
Let $x_1(t)$ be a smooth real-valued function of $t$ and let
$z_1(t)$, $z_2(t)$, $z_3(t)$ be smooth complex-valued functions of
$t$ such that
\begin{align}
\label{ch4s2eq9}
\frac{dx_1}{dt} & =  0\,\text{,} \\[4pt]
\label{ch4s2eq10}
\frac{dz_1}{dt} & = -\nu z_1 -\lambda \overline{z_2z_3}\,\text{,} \\[4pt]
\label{ch4s2eq11}
\frac{dz_2}{dt} & =  \mu z_2 -\lambda \overline{z_3z_1}\;\text{and} \\[4pt]
\label{ch4s2eq12} \frac{dz_3}{dt} & =  (\nu - \mu) z_3 -\lambda
\overline{z_1z_2}\,\text{,}
\end{align}
using the notation from Definition \ref{ch4s2dfn1}.  There exists
$\epsilon
> 0$ such that these equations have a solution for
$t\in(-\epsilon,\epsilon)$ and the subset $M$ of\/
$\R\oplus\C^3\cong\R^7$ defined by
\begin{align*}
M  = \big\{\big(x_1(t)+c,\,e^{i\phi_1}z_1(t),\,e^{i\phi_2}z_2(t),
\,e^{-i(\phi_1+\phi_2)}&z_3(t)\big):\\
&\,t\in(-\epsilon,\epsilon),\,(c,e^{i\phi_1}, e^{i\phi_2})\in
\text{\emph{G}}\big\}
\end{align*}
is an associative 3-fold in $\R^7$.  Moreover, $M$ does not lie in
$\{x\}\times\C^3$ for any\/ $x\in\R$, as long as not both $\mu$ and
$\nu$ are zero, and \eq{ch4s2eq10}-\eq{ch4s2eq12} imply that
$\,\Im(z_1z_2z_3)=A$, where $A$ is a real constant.
\end{thm}

\begin{proof} We only need to prove the last sentence in the statement above.
We deduce immediately from \eq{ch4s2eq9} that $x_1$ is constant in
the direction transverse to the group action, though it \emph{is}
changing along the group action (as long as not both $\mu$ and $\nu$
are zero), which means that $M$ does not lie in $\{ x \} \times
\C^3$ for any real constant $x$ in this case. We also note from
\eq{ch4s2eq10}-\eq{ch4s2eq12} that
\begin{equation*}
\frac{d}{dt}\,(z_1z_2z_3)  =
-\lambda(|z_2|^2|z_3|^2+|z_3|^2|z_1|^2+|z_1|^2|z_2|^2),
\end{equation*}
which is real, therefore $\Im(z_1z_2z_3)$ is a real constant.
\end{proof}

There are two trivial cases which may be solved immediately.

Firstly, suppose $\lambda=0$.  This is not geometrically interesting
since it implies that $\text{G}$ contains all possible translations
in the first coordinate.  Solving \eq{ch4s2eq10}-\eq{ch4s2eq12}
shows that
\begin{align*}
M= \R \times \Big\{\Big(A_1e^{i\phi_1 -\nu t}, \, A_2e^{i\phi_2 +
\mu t}, \, A_3&e^{-i(\phi_1 + \phi_2) + (\nu - \mu) t}\Big) :\,\\
&\qquad\qquad\quad t \in \R, \, \mu \phi_1 + \nu \phi_2 = 0 \Big\},
\end{align*}
\noindent where $A_1$, $A_2$, $A_3$ are complex constants such that
$\Im(A_1A_2 A_3)=A$.  The expression in brackets above defines a
holomorphic curve in $\C^3$.

The other case is when $\mu=\nu=0$.  This forces $c = 0$ in
$\text{G}$, so there is no translation action in $\text{G}$, which
means that $M$ will be an embedded $\U(1)^2$-invariant SL 3-fold as
studied in \cite[$\S$III.3.A]{HarLaw}:
\begin{align*}
M=\{(x_1,z_1,z_2,z_3)\in\R^7:\,x_1=x,\, \Im(z_1&z_2z_3)=A, \,\\
&|z_1|^2 - |z_3|^2 = B, \, |z_2|^2 - |z_3|^2 = C\}
\end{align*}
\noindent for some $x,A,B,C\in\R$.

\subsection{$\U(1)^2$-invariant coassociative cones}

We consider coassociative 4-folds invariant both under the action of
$\U(1)^2$ on the $\C^3$ component of $\R^7\cong\R\oplus\C^3$ and
under dilations.

\begin{dfn}\label{ch5s2subs1dfn1}
Let $\R^+$ denote the group of positive real numbers under
multiplication.   Define an action of $\R^+\times\U(1)^2$ on
$\R^7\cong\R\oplus\C^3$ by
\begin{equation}\label{ch5s2subs1eq1}
(x_1,z_1,z_2,z_3)\longmapsto(rx_1,\,re^{i\phi_1}z_1,\,re^{i\phi_2}z_2,\,re^{-i(\phi_1+\phi_2)}z_3),
\quad r>0, \, \phi_1,\phi_2\in\R.
\end{equation}
\end{dfn}

\vspace{-12pt}

We again apply the method described in $\S$\ref{symmethod}, though
this time we must choose our orbit so that $\varphi_0$ vanishes on
it. This constraint imposes the condition $\Re(z_1z_2z_3)=0$.  We
thus have the following result.

\begin{thm}\label{ch5s2subs1thm1}
Let $x_1(t)$ be a smooth real-valued function of $t$ and let
$z_1(t)$, $z_2(t)$, $z_3(t)$ be smooth complex-valued functions of
$t$ satisfying
\begin{align}
\label{ch5s2subs1eq6}
\frac{dx_1}{dt} & =  -3{\rm Im}(z_1z_2z_3), \\[4pt]
\label{ch5s2subs1eq7}
\frac{dz_1}{dt} & =  z_1(|z_2|^2-|z_3|^2)+ix_1\overline{z_2z_3},
\\[4pt]
\label{ch5s2subs1eq8}
\frac{dz_2}{dt} & =
z_2(|z_3|^2-|z_1|^2)+ix_1\overline{z_3z_1}\,\text{and}
\\[4pt]
\label{ch5s2subs1eq9}
\frac{dz_3}{dt} & = z_3(|z_1|^2-|z_2|^2)+ix_1\overline{z_1z_2},
\end{align}
\noindent along with the condition
\begin{equation}\label{ch5s2subs1eq10}
\Re(z_1z_2z_3) = 0
\end{equation}
at\/ $t=0$.  The subset $M$ of\/ $\R\oplus\C^3\cong\R^7$ defined by
\begin{equation*}
M = \left\{ \big(rx_1(t),\, re^{i\phi_1}z_1(t),\, re^{i\phi_2}z_2(t),\,
re^{-i(\phi_1+\phi_2)}z_3(t)\big)\,:\, r>0, \, \phi_1,
\phi_2,t\in\R\right\}
\end{equation*}
\noindent is a coassociative 4-fold in $\R^7$.  Moreover,
\eq{ch5s2subs1eq10} holds for all $t\in\R$ and
$x_1^2+|z_1|^2+|z_2|^2+|z_3|^2$ is a constant which can be taken to
be $1$.
\end{thm}

\begin{proof}
It is immediate from \eq{ch5s2subs1eq6}-\eq{ch5s2subs1eq9} that
$x_1^2+|z_1|^2+|z_2|^2+|z_3|^2$ is a constant which can be chosen to
be $1$ without loss of generality.  We may also calculate
$$\frac{d}{dt}(z_1z_2z_3)=ix_1(|z_2|^2|z_3|^2+|z_3|^2|z_1|^2+|z_1|^2|z_2|^2)$$
using \eq{ch5s2subs1eq6}-\eq{ch5s2subs1eq9} and deduce that
$\Re(z_1z_2z_3)$ is a constant which has to be zero since
\eq{ch5s2subs1eq10} holds at $t=0$.  Theorem \ref{ch5s1thm4} only
gives us that solutions to \eq{ch5s2subs1eq6}-\eq{ch5s2subs1eq9}
exist for $t\in(-\epsilon,\epsilon)$ for some $\epsilon>0$, but
solutions exist for all $t$, as argued in the proof of Theorem
\ref{ch4s2thm2}, since the functions involved are all bounded.
\end{proof}

\subsection{$\U(1)^2$-invariant Cayley cones}

We conclude by turning our attention to Cayley cones which are
invariant under a $\U(1)^2$ subgroup of $\U(1)^4$.

\begin{dfn}\label{ch5s3subs1dfn1}
Let $\text{G}\subseteq\U(1)^4$ be defined by
\begin{align*}
\text{G}  =  \big\{
(e^{i\alpha_1}&,\,e^{i\alpha_2},\,e^{i\alpha_3},\,e^{i\alpha_4})\,:\,
\alpha_1,\alpha_2,\,\alpha_3,\alpha_4\in\R\,\;\text{satisfy}\\
&\qquad\alpha_1+\alpha_2+\alpha_3+\alpha_4=0\,\;\text{and}\,\;
a_1\alpha_1+a_2\alpha_2+ a_3\alpha_3+a_4\alpha_4=0\big\}
\end{align*}
for coprime integers $a_1,a_2,a_3,a_4$ with $a_1+a_2+a_3+a_4=0$ and
$a_1\leq a_2\leq a_3\leq a_4$.  This acts on $\C^4\cong\R^8$ in the
obvious way as a $\U(1)^2$ subgroup of $\U(1)^4$.
\end{dfn}

\begin{thm} Use the notation of Definition \ref{ch5s3subs1dfn1}.
Let $z_j(t)$ for $j=1,2,3,4$ be smooth complex-valued functions of
$t$ satisfying
\begin{align}
\label{ch5s3subs1eq4}
\frac{dz_1}{dt} & =
a_1\overline{z_2z_3z_4}+\frac{1}{2}\,z_1\big((a_4-a_3)|z_2|^2
+(a_2-a_4)|z_3|^2+(a_3-a_2)|z_4|^2\big), \\[4pt]
\label{ch5s3subs1eq5}
\frac{dz_2}{dt} & =
a_2\overline{z_3z_4z_1}+\frac{1}{2}\,z_2\big((a_4-a_1)|z_3|^2
+(a_1-a_3)|z_4|^2+(a_3-a_4)|z_1|^2\big), \displaybreak[0]\\[4pt]
\label{ch5s3subs1eq6}
\frac{dz_3}{dt} & =
a_3\overline{z_4z_1z_2}+\frac{1}{2}\,z_3\big((a_2-a_1)|z_4|^2
+(a_4-a_2)|z_1|^2+(a_1-a_4)|z_2|^2\big)\,\text{and} \\[4pt]
\label{ch5s3subs1eq7}
\frac{dz_4}{dt} & =
a_4\overline{z_1z_2z_3}+\frac{1}{2}\,z_4\big((a_2-a_3)|z_1|^2
+(a_3-a_1)|z_2|^2+(a_1-a_2)|z_3|^2\big).
\end{align}
The subset $M$ of\/ $\C^4\cong\R^8$ given by
\begin{align*}
M=\big\{\big(re^{i\alpha_1}z_1(t),\,re^{i\alpha_2}z_2(t),\,re^{i\alpha_3}z_3&(t),\,re^{i\alpha_4}z_4(t)\big)\,:\,\\
&r>0,\,(e^{i\alpha_1},e^{i\alpha_2},e^{i\alpha_3},e^{i\alpha_4})\in\text{\emph{G}},\,t\in\R\big\}
\end{align*}
is a Cayley 4-fold in $\R^8$.  Moreover,
$|z_1|^2+|z_2|^2+|z_3|^2+|z_4|^2$ is a constant which can be taken
to be $1$ and $\Im(z_1z_2z_3z_4)=A$ for some real constant $A$.
\end{thm}

\begin{proof}  It is clear from
\eq{ch5s3subs1eq4}-\eq{ch5s3subs1eq7} that $|z_1|^2+\ldots+|z_4|^2$
is a constant and that we can take this constant to be $1$ without
loss of generality.  Furthermore,
\begin{equation*}
\frac{d}{dt}\,(z_1z_2z_3z_4)=a_1|z_2z_3z_4|^2+a_2|z_3z_4z_1|^2+a_3|z_4z_1z_2|^2+
a_4|z_1z_2z_3|^2,
\end{equation*}
which is purely real.  Therefore $\Im(z_1z_2z_3z_4)=A$ is constant.
Theorem \ref{ch5s1thm5} only gives existence of solutions of
$t\in(-\epsilon,\epsilon)$ for some $\epsilon>0$. However, by the
same argument as in the proof of Theorem \ref{ch4s2thm2}, solutions
exist for all $t\in\R$, using the boundedness of the functions
involved.
\end{proof}

\begin{ack}
I am indebted to Dominic Joyce for his help with this research
project.  I would also like to thank Alexei Kovalev and Andrew
Dancer for useful corrections and suggestions and EPSRC for
providing the funding for this study.  
\end{ack}


\end{document}